# IMPROVED MINIMAX PREDICTIVE DENSITIES UNDER KULLBACK–LEIBLER LOSS[1]

By Edward I. George, Feng Liang and Xinyi Xu

*University of Pennsylvania, Duke University
and Ohio State University*

Let $X|\mu \sim N_p(\mu, v_x I)$ and $Y|\mu \sim N_p(\mu, v_y I)$ be independent $p$-dimensional multivariate normal vectors with common unknown mean $\mu$. Based on only observing $X = x$, we consider the problem of obtaining a predictive density $\hat{p}(y|x)$ for $Y$ that is close to $p(y|\mu)$ as measured by expected Kullback–Leibler loss. A natural procedure for this problem is the (formal) Bayes predictive density $\hat{p}_U(y|x)$ under the uniform prior $\pi_U(\mu) \equiv 1$, which is best invariant and minimax. We show that any Bayes predictive density will be minimax if it is obtained by a prior yielding a marginal that is superharmonic or whose square root is superharmonic. This yields wide classes of minimax procedures that dominate $\hat{p}_U(y|x)$, including Bayes predictive densities under superharmonic priors. Fundamental similarities and differences with the parallel theory of estimating a multivariate normal mean under quadratic loss are described.

**1. Introduction.** Let $X|\mu \sim N_p(\mu, v_x I)$ and $Y|\mu \sim N_p(\mu, v_y I)$ be independent $p$-dimensional multivariate normal vectors with common unknown mean $\mu$, and let $p(x|\mu)$ and $p(y|\mu)$ denote the conditional densities of $X$ and $Y$. We assume that $v_x$ and $v_y$ are known.

Based on only observing $X = x$, we consider the problem of obtaining a predictive density $\hat{p}(y|x)$ for $Y$ that is close to $p(y|\mu)$. We measure this closeness by Kullback–Leibler (KL) loss,

$$L(\mu, \hat{p}(\cdot|x)) = \int p(y|\mu) \log \frac{p(y|\mu)}{\hat{p}(y|x)} \, dy, \tag{1}$$

Received July 2003; revised March 2005.
[1]Supported by NSF Grant DMS-01-30819.

and evaluate $\hat{p}$ by its expected loss or risk function

$$R_{\text{KL}}(\mu, \hat{p}) = \int p(x|\mu) L(\mu, \hat{p}(\cdot|x))\, dx. \tag{2}$$

For the comparison of two procedures, we say that $\hat{p}_1$ dominates $\hat{p}_2$ if $R_{\text{KL}}(\mu, \hat{p}_1) \leq R_{\text{KL}}(\mu, \hat{p}_2)$ for all $\mu$ and with strict inequality for some $\mu$. By a sufficiency and transformation reduction, this problem is seen to be equivalent to estimating the predictive density of $X_{n+1}$ under KL loss based on observing $X_1, \ldots, X_n$ when $X_1, \ldots, X_{n+1}|\mu$ i.i.d. $\sim N_p(\mu, \Sigma)$. For distributions beyond the normal, versions and approaches for the KL risk prediction problem have been developed by Aslan [2], Harris [10], Hartigan [11], Komaki [12, 14] and Sweeting, Datta and Ghosh [24].

For any prior distribution $\pi$ on $\mu$, Aitchison [1] showed that the average risk $r(\pi, \hat{p}) = \int R_{\text{KL}}(\mu, \hat{p}) \pi(\mu)\, d\mu$ is minimized by

$$\hat{p}_\pi(y|x) = \int p(y|\mu) \pi(\mu|x)\, d\mu, \tag{3}$$

which we will refer to as a Bayes predictive density. Unless $\pi$ is a trivial point prior, $\hat{p}_\pi(y|x) \notin \{p(y|\mu) : \mu \in R^p\}$, that is, $\hat{p}_\pi$ will not correspond to a "plug-in" estimate for $\mu$, although under suitable conditions on $\pi$, $\hat{p}_\pi(y|x) \to p(y|\mu)$ as $v_x \to 0$.

For this problem, the best invariant predictive density (with respect to the location group) is the Bayes predictive density under the uniform prior $\pi_{\text{U}}(\mu) \equiv 1$, namely

$$\hat{p}_{\text{U}}(y|x) = \frac{1}{\{2\pi(v_x + v_y)\}^{p/2}} \exp\left\{-\frac{\|y - x\|^2}{2(v_x + v_y)}\right\}, \tag{4}$$

which has constant risk; see [18] and [19]. More precisely, one might refer to $\hat{p}_{\text{U}}$ as a formal Bayes procedure because $\pi_{\text{U}}$ is improper. Aitchison [1] showed that $\hat{p}_{\text{U}}(y|x)$ dominates the plug-in predictive density $p(y|\hat{\mu}_{\text{MLE}})$ which simply substitutes the maximum likelihood estimate $\hat{\mu}_{\text{MLE}} = x$ for $\mu$. As will be seen in Section 2, $\hat{p}_{\text{U}}$ is minimax for KL loss (1). That $\hat{p}_{\text{U}}$ is best invariant and minimax can also be seen as a special case of the more general recent results in Liang and Barron [17], who also show that $\hat{p}_{\text{U}}$ is admissible when $p = 1$ under the same loss.

However, $\hat{p}_{\text{U}}$ is inadmissible when $p \geq 3$. Komaki [13] proved that when $p \geq 3$, $\hat{p}_{\text{U}}$ itself is dominated by the (formal) Bayes predictive density

$$\hat{p}_{\text{H}}(y|x) = \int p(y|\mu) \pi_{\text{H}}(\mu|x)\, d\mu, \tag{5}$$

where

$$\pi_{\text{H}}(\mu) = \|\mu\|^{-(p-2)} \tag{6}$$

is the (improper) harmonic prior recommended by Stein [21], which we subscript by "H" for harmonic. Although Komaki referred to $\pi_H$ as harmonic, his proof did not directly exploit this property.

More recently, Liang [16] showed that $\hat{p}_U$ is also dominated by the proper Bayes predictive density $\hat{p}_a(y|x)$ under the prior $\pi_a(\mu)$ (see [23]) defined hierarchically as

$$(7) \qquad \mu|s \sim N_p(0, sv_0 I), \qquad s \sim (1+s)^{a-2}.$$

Here $v_0$ and $a$ are hyperparameters. The conditions for domination are that $v_0 \geq v_x$, and $a \in [0.5, 1)$ when $p = 5$ and $a \in [0, 1)$ when $p \geq 6$. Note that $\pi_a$ depends on the constant $v_0$ in (7), a dependence that will be maintained throughout this paper. The harmonic prior $\pi_H$ is well known to be the special case of $\pi_a$ when $a = 2$.

These results closely parallel some key developments concerning minimax estimation of a multivariate normal mean under quadratic loss. Based on observing $X|\mu \sim N_p(\mu, I)$, that problem is to estimate $\mu$ under

$$(8) \qquad R_Q(\mu, \hat{\mu}) = E_\mu \|\hat{\mu} - \mu\|^2,$$

where we have denoted quadratic risk by $R_Q$ to distinguish it from the KL risk $R_{KL}$ in (2). Under $R_Q$, $\hat{\mu}_{MLE} = X$ is best invariant and minimax, and is admissible if and only if $p \leq 2$. Note that $\hat{\mu}_{MLE}$ plays the same role here that $\hat{p}_U$ plays in our KL risk problem. A further connection between $\hat{\mu}_{MLE}$ and $\hat{p}_U$ is revealed by the fact that $\hat{\mu}_{MLE} \equiv E_{\pi_U}(\mu|x)$, the posterior mean of $\mu$ under $\pi_U(\mu) \equiv 1$.

Stein [21] showed that $\hat{\mu}_H = E_{\pi_H}(\mu|x)$, the posterior mean under $\pi_H$, dominates $\hat{\mu}_{MLE}$ when $p \geq 3$, and Strawderman [23] showed that $\hat{\mu}_a = E_{\pi_a}(\mu|x)$, the proper Bayes rule under $\pi_a$ when $v_x = v_0 = 1$, dominates $\hat{\mu}_{MLE}$ when $a \in [0.5, 1)$ for $p = 5$ and when $a \in [0, 1)$ for $p \geq 6$. Comparing these results to those of Komaki and Liang in the predictive density problem, the parallels are striking. A principal purpose of our paper is to draw out these parallels in a more unified and transparent way.

For these and other shrinkage domination results in the quadratic risk estimation problem, there exists a unifying theory that focuses on the properties of the marginal distribution of $X$ under $\pi$, namely

$$(9) \qquad m_\pi(x) = \int p(x|\mu) \pi(\mu) \, d\mu.$$

The key to this theory is the representation due to Brown [4] that any posterior mean of $\mu$, $\hat{\mu}_\pi = E_\pi(\mu|x)$, is of the form

$$(10) \qquad \hat{\mu}_\pi = x + \nabla \log m_\pi(x),$$

where $\nabla = (\partial/\partial x_1, \ldots, \partial/\partial x_p)'$. To show that $\hat{\mu}_H$ dominates $\hat{\mu}_{MLE}$, Stein [21, 22] used this representation to establish that $R_Q(\mu, \hat{\mu}_{MLE}) - R_Q(\mu, \hat{\mu}_\pi) = E_\mu U(X)$, where

$$U(X) = \|\nabla \log m_\pi(X)\|^2 - 2\frac{\nabla^2 m_\pi(X)}{m_\pi(X)} \tag{11}$$

$$= -4\frac{\nabla^2 \sqrt{m_\pi(X)}}{\sqrt{m_\pi(X)}} \tag{12}$$

is an unbiased estimate of the risk reduction of $\hat{\mu}_\pi$ over $\hat{\mu}_{MLE}$, where $\nabla^2 m_\pi(x) = \sum \frac{\partial^2}{\partial x_i^2} m_\pi(x)$.

Because $\hat{\mu}_{MLE}$ is minimax, it follows immediately from (11) that $\nabla^2 m_\pi(x) \leq 0$ is a sufficient (though not necessary) condition for $\hat{\mu}_\pi$ to be minimax, and as long as $m_\pi(x)$ is not constant, for $\hat{\mu}_\pi$ to dominate $\hat{\mu}_{MLE}$. [Recall that a function $m(x)$ is superharmonic when $\nabla^2 m(x) \leq 0$.] The fact that $\hat{\mu}_H$ dominates $\hat{\mu}_{MLE}$ when $p \geq 3$ now follows easily from the fact that nonconstant superharmonic priors [of which the harmonic prior $\pi_H(\mu)$ is of course a special case] yield superharmonic marginals $m_\pi$ for $X$.

It follows from (12) that the weaker condition $\nabla^2 \sqrt{m_\pi(x)} \leq 0$ is sufficient for $\hat{\mu}_\pi$ to be minimax, although strict inequality on a set of positive Lebesgue measure is then needed to guarantee domination over $\hat{\mu}_{MLE}$. Fourdrinier, Strawderman and Wells [6] showed that the Strawderman priors $\pi_a$ in (7) yield superharmonic $\sqrt{m_\pi}$, so that the minimaxity of the Strawderman estimators is established by (12). In fact, it follows from their results that $\pi_a$ also yields superharmonic $\sqrt{m_\pi}$ when $a \in [1, 2)$ and $p \geq 3$, thereby broadening the class of formal Bayes minimax estimators.

One major aim of the present paper is to establish an analogous unifying theory for the KL risk prediction problem. Paralleling (10), we begin by showing how any Bayes predictive density $\hat{p}_\pi$ can be explicitly represented in terms of $\hat{p}_U$ and the form of the corresponding marginal $m_\pi$. Coupled with the heat equation, Brown's representation and Stein's identity, this representation is seen to lead to a new identity that links KL risk reduction to Stein's unbiased estimate of risk reduction. Based on this link, we obtain sufficient conditions on $m_\pi$ for minimaxity and domination of $\hat{p}_\pi$ over $\hat{p}_U$. These general conditions subsume the specialized results of Komaki [13] and Liang [16] and can be used to obtain wide classes of improved minimax Bayes predictive densities including $\hat{p}_H$ and $\hat{p}_a$. Furthermore, the underlying priors and marginals can be readily adapted to obtain minimax shrinkage toward an arbitrary point or subspace, and linear combinations of superharmonic priors and marginals can be constructed to obtain minimax multiple shrinkage predictive density analogues of the minimax multiple shrinkage estimators of George [7, 8, 9]. Thus, the parallels between the estimation and the prediction problem are broad, both qualitatively and technically. The main contribution of this paper is to establish this interesting connection.

**2. General conditions for minimaxity.** In this section we develop and prove our main results concerning general conditions under which a Bayes predictive density $\hat{p}_\pi(y|x)$ in (3) will be minimax and dominate $\hat{p}_U(y|x)$. We begin with three lemmas that may also be of independent interest. The following general notation will be useful throughout. For $Z|\mu \sim N_p(\mu, vI)$ and a prior $\pi$ on $\mu$, we denote the marginal distribution of $Z$ by

$$m_\pi(z; v) = \int p(z|\mu)\pi(\mu)\,d\mu. \tag{13}$$

In terms of this notation, the marginal distributions of $X|\mu \sim N_p(\mu, v_x I)$ and $Y|\mu \sim N_p(\mu, v_y I)$ under $\pi$ are then $m_\pi(x; v_x)$ and $m_\pi(y; v_y)$, respectively.

LEMMA 1. *If $m_\pi(z; v_x)$ is finite for all $z$, then for every $x$, $\hat{p}_\pi(y|x)$ will be a proper probability distribution over $y$. Furthermore, the mean of $\hat{p}_\pi(y|x)$ is equal to $E_\pi(\mu|x)$.*

PROOF. Both claims follow by integrating (3) with respect to $y$ and switching the order of integration using the Fubini–Tonelli theorem. □

Lemma 1 is important because, for our decision problem to be meaningful, it is necessary for a predictive density to be a proper probability distribution. By the laws of probability, a Bayes predictive density $\hat{p}_\pi(y|x)$ will be a proper probability distribution whenever $\pi(\mu)$ is a proper prior distribution. But by Lemma 1, improper $\pi(\mu)$ can still yield proper $\hat{p}_\pi(y|x)$ under a very weak condition.

Our next lemma establishes a key alternative representation of $\hat{p}_\pi(y|x)$ that makes use of the weighted mean

$$W = \frac{v_y X + v_x Y}{v_x + v_y}. \tag{14}$$

Note that $W$ would be a sufficient statistic for $\mu$ if both $X$ and $Y$ were observed. As $X$ and $Y$ are independent (conditionally on $\mu$), it follows that $W|\mu \sim N_p(\mu, v_w I)$ where

$$v_w = \frac{v_x v_y}{v_x + v_y}.$$

The marginal distribution of $W$ is then $m_\pi(w; v_w)$.

LEMMA 2. *For any prior $\pi(\mu)$, $\hat{p}_\pi(y|x)$ can be expressed as*

$$\hat{p}_\pi(y|x) = \frac{m_\pi(w; v_w)}{m_\pi(x; v_x)} \hat{p}_U(y|x), \tag{15}$$

where $\hat{p}_{\mathrm{U}}(y|x)$ is defined by (4). Furthermore, the difference between the KL risks of $\hat{p}_{\mathrm{U}}(y|x)$ and $\hat{p}_\pi(y|x)$ is given by

$$
\begin{aligned}
R_{\mathrm{KL}}(\mu, \hat{p}_{\mathrm{U}}) &- R_{\mathrm{KL}}(\mu, \hat{p}_\pi) \\
&= E_{\mu, v_w} \log m_\pi(W; v_w) - E_{\mu, v_x} \log m_\pi(X; v_x),
\end{aligned}
\tag{16}
$$

where $E_{\mu,v}(\cdot)$ stands for expectation with respect to the $N(\mu, vI)$ distribution.

PROOF. The joint marginal distribution of $X$ and $Y$ under $\pi$ is,

$$
\begin{aligned}
p_\pi(x, y) &= \int p(x|\mu) p(y|\mu) \pi(\mu) \, d\mu \\
&= \int \frac{1}{(2\pi v_x)^{p/2}} \exp\left\{-\frac{\|x-\mu\|^2}{2v_x}\right\} \\
&\quad \times \frac{1}{(2\pi v_y)^{p/2}} \exp\left\{-\frac{\|y-\mu\|^2}{2v_y}\right\} \pi(\mu) \, d\mu \\
&= \int \frac{1}{\{2\pi(v_x + v_y)\}^{p/2}} \exp\left\{-\frac{\|y-x\|^2}{2(v_x + v_y)}\right\} \\
&\quad \times \frac{1}{(2\pi v_w)^{p/2}} \exp\left\{-\frac{\|w-\mu\|^2}{2v_w}\right\} \pi(\mu) \, d\mu \\
&= \hat{p}_{\mathrm{U}}(y|x) m_\pi(w; v_w).
\end{aligned}
$$

The representation (15) now follows since $\hat{p}_\pi(y|x) = p_\pi(x,y)/m_\pi(x; v_x)$.

To prove (16), the KL risk difference can be expressed as

$$
\begin{aligned}
R_{\mathrm{KL}}(\mu, \hat{p}_{\mathrm{U}}) - R_{\mathrm{KL}}(\mu, \hat{p}_\pi) &= \int \int p(x|\mu) p(y|\mu) \log \frac{\hat{p}_\pi(y|x)}{\hat{p}_{\mathrm{U}}(y|x)} \, dx \, dy \\
&= \int \int p(x|\mu) p(y|\mu) \log \frac{m_\pi(w; v_w)}{m_\pi(x; v_x)} \, dx \, dy,
\end{aligned}
$$

where the second equality makes use of (15). The second expression in (16) is seen to equal this last expression by the change of variable theorem. □

Paralleling Brown's representation (10), representation (15) reveals the explicit role played by the marginal distribution of the data under $\pi$. Analogous to Bayes estimators $E_\pi(\mu|x)$ of $\mu$ that "shrink" $\hat{\mu}_{\mathrm{MLE}} = x$, this representation reveals that Bayes predictive densities $\hat{p}_\pi(y|x)$ "shrink" $\hat{p}_{\mathrm{U}}(y|x)$ by a factor $m_\pi(w; v_w)/m_\pi(x; v_x)$. However, the nature of the shrinkage by $\hat{p}_\pi(y|x)$ is different than that by $E_\pi(\mu|x)$. To insure that $\hat{p}_\pi(y|x)$ remains a proper probability distribution, the factor $m_\pi(w; v_w)/m_\pi(x; v_x)$ cannot be strictly less than 1. In contrast to simply shifting $\hat{\mu}_{\mathrm{MLE}} = x$ toward the mean of $\pi$, $\hat{p}_\pi(y|x)$ adjusts $\hat{p}_{\mathrm{U}}(y|x)$ to concentrate more on the higher probability regions of $\pi$. Figure 1 illustrates such shrinkage of $\hat{p}_{\mathrm{U}}(y|x)$ by $\hat{p}_{\mathrm{H}}(y|x)$ in (5) when $v_x = 1$, $v_y = 0.2$ and $p = 5$.

For our purposes, the principal benefit of (15) is that it reduces the KL risk difference (16) to a simple functional of the marginal $m_\pi(z;v)$. As will be seen in the proof of Theorem 1 below, (16) is the key to establishing general conditions for the dominance of $\hat{p}_\pi$ over $\hat{p}_{\mathrm{U}}$. First, however, we use it to facilitate a simple direct proof of the minimaxity of $\hat{p}_{\mathrm{U}}$, a result that also follows from the more general results of Liang and Barron [17].

COROLLARY 1. *The Bayes predictive density under $\pi(\mu) \equiv 1$, namely $\hat{p}_{\mathrm{U}}$, is minimax under $R_{\mathrm{KL}}$.*

PROOF. By a transformation of variables, $x \to (x-\mu)$ and $y \to (y-\mu)$, it is easy to see that $R_{\mathrm{KL}}(\mu, \hat{p}_{\mathrm{U}}) = R_{\mathrm{KL}}(0, \hat{p}_{\mathrm{U}}) = r$ for all $\mu$, so that $R_{\mathrm{KL}}(\mu, \hat{p}_{\mathrm{U}})$ is constant. Next, we show that $r$ is a Bayes risk limit of a sequence of Bayes rules $\hat{p}_{\pi_n}$ with $\pi_n(\mu) = N_p(0, \sigma_n^2 I)$, where $\sigma_n^2 \to \infty$ as $n \to \infty$. By the fact that $r(\pi_n, \hat{p}_{\mathrm{U}}) \equiv r$ and (16),

$$(17) \quad r - r(\pi_n, \hat{p}_{\pi_n}) = \int \pi_n(\mu) [E_{\mu, v_w} \log m_{\pi_n}(W; v_w) - E_{\mu, v_x} \log m_{\pi_n}(X; v_x)] \, d\mu,$$

where

$$m_{\pi_n}(z; v) = (2\pi(v + \sigma_n^2))^{-p/2} \exp\left\{-\frac{\|z\|^2}{2(v + \sigma_n^2)}\right\}.$$

It is now easy to check that $(17) = O(1/\sigma_n^2)$ and hence goes to zero as $n$ goes to infinity. By Theorem 5.18 of [3], the minimaxity of $\hat{p}_{\mathrm{U}}$ follows. □

Our next lemma provides a new identity that links $E_{\mu,v} \log m_\pi(Z;v)$ to Stein's unbiased estimate of risk reduction $U(x)$ in (11) and (12) for the quadratic risk estimation problem. When combined with (16) in Theorem 1, this identity will be seen to play a key role in establishing sufficient conditions on $m_\pi$ for $\hat{p}_\pi$ to be minimax and to dominate $\hat{p}_{\mathrm{U}}$.

LEMMA 3. *If $m_\pi(z; v_x)$ is finite for all $z$, then for any $v_w \leq v \leq v_x$, $m_\pi(z; v)$ is finite. Moreover,*

$$(18) \quad \frac{\partial}{\partial v} E_{\mu,v} \log m_\pi(Z;v) = E_{\mu,v}\left(\frac{\nabla^2 m_\pi(Z;v)}{m_\pi(Z;v)} - \frac{1}{2}\|\nabla \log m_\pi(Z;v)\|^2\right)$$

$$(19) \quad = E_{\mu,v}\left(2\frac{\nabla^2 \sqrt{m_\pi(Z;v)}}{\sqrt{m_\pi(Z;v)}}\right).$$

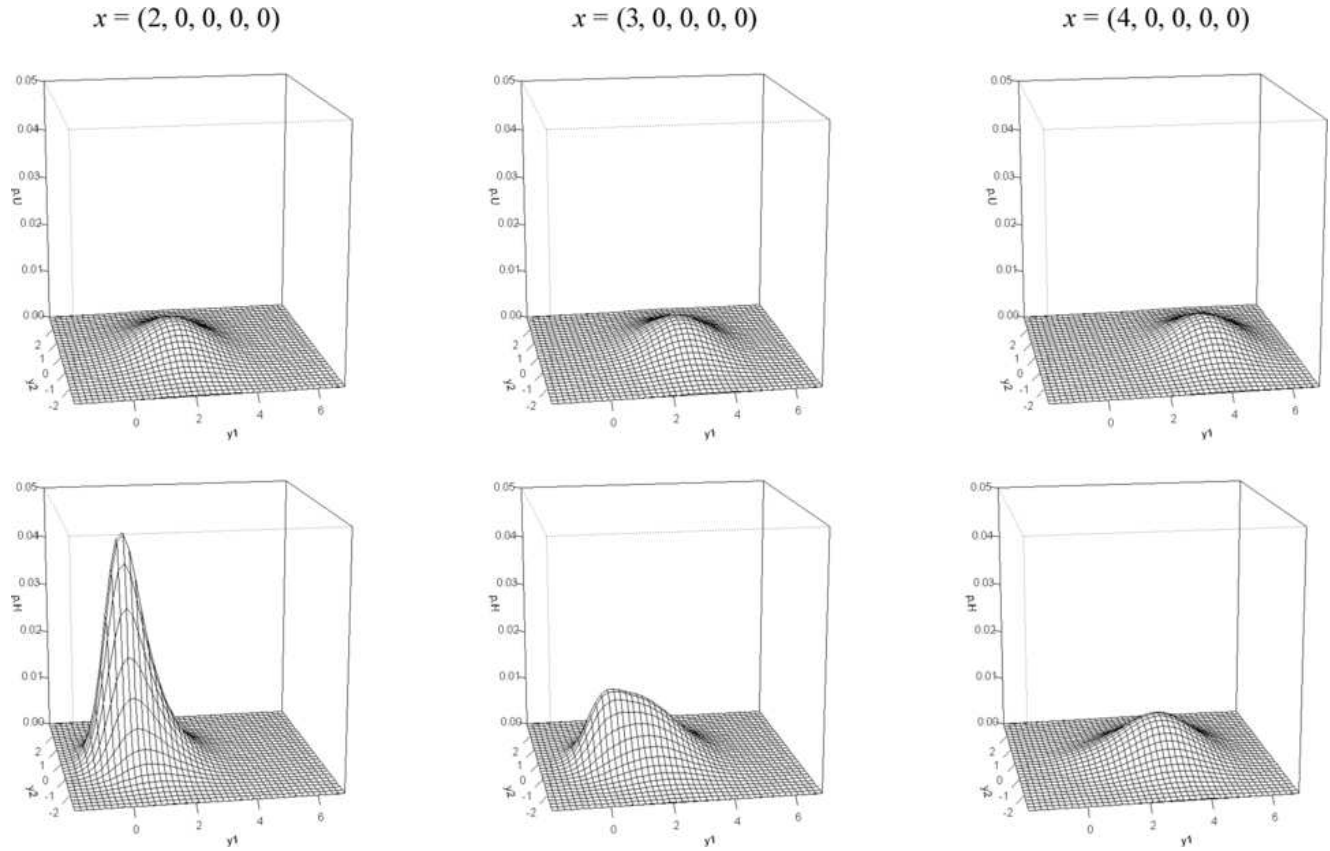

FIG. 1. *Shrinkage of $\hat{p}_U(y|x)$ to obtain $\hat{p}_H(y|x)$ when $v_x = 1, v_y = 0.2$ and $p = 5$. Here $y = (y_1, y_2, 0, 0, 0)$.*

PROOF. When $m_\pi(z; v_x)$ is finite for all $z$, it is easy to check that for any fixed $z$ and any $v_w \leq v \leq v_x$,

$$m_\pi(z; v) \leq \left(\frac{v_x}{v_w}\right)^{p/2} m_\pi(z; v_x) < \infty.$$

Letting $Z^* = (Z - \mu)/\sqrt{v} \sim N_p(0, I)$, we obtain

(20)
$$\frac{\partial}{\partial v} E_{\mu,v} \log m_\pi(Z; v) = \frac{\partial}{\partial v} E \log m_\pi(\sqrt{v} Z^* + \mu; v)$$
$$= E \frac{(\partial/\partial v) m_\pi(\sqrt{v} Z^* + \mu; v)}{m_\pi(\sqrt{v} Z^* + \mu; v)},$$

where

$$\frac{\partial}{\partial v} m_\pi(\sqrt{v} z^* + \mu; v)$$
$$= \frac{\partial}{\partial v} \int \frac{1}{(2\pi v)^{p/2}} \exp\left\{-\frac{\|\sqrt{v} z^* + \mu - \mu'\|^2}{2v}\right\} \pi(\mu') \, d\mu'$$
$$= \int \left(-\frac{p}{2v} + \frac{\|z - \mu'\|^2}{2v^2} - \frac{\|z^*\|^2}{2v} - \frac{z^{*\prime}(\mu - \mu')}{2v^{3/2}}\right) p(z|\mu') \pi(\mu') \, d\mu'$$
$$= \frac{\partial}{\partial v} m_\pi(z; v) - \int \frac{(z - \mu)'(z - \mu')}{2v^2} p(z|\mu') \pi(\mu') \, d\mu'.$$

Using the fact that

(21)
$$\frac{\partial}{\partial v} m_\pi(z; v) = \frac{1}{2} \nabla^2 m_\pi(z; v),$$

which is straightforward to verify, and by Brown's representation $E_\pi(\mu'|z) = z + v \nabla \log m_\pi(z)$ from (10),

(22)
$$E \frac{(\partial/\partial v) m_\pi(\sqrt{v} Z^* + \mu; v)}{m_\pi(\sqrt{v} Z^* + \mu; v)}$$
$$= E_{\mu,v} \left(\frac{1}{2} \frac{\nabla^2 m_\pi(Z; v)}{m_\pi(Z; v)} + \frac{(Z - \mu)' \nabla \log m_\pi(Z; v)}{2v}\right).$$

Finally, by (2.3) of [22],

(23)
$$E_{\mu,v} \frac{(Z - \mu)' \nabla \log m_\pi(Z; v)}{2v}$$
$$= E_{\mu,v} \frac{1}{2} \nabla^2 \log m_\pi(Z; v) = E_{\mu,v} \frac{1}{2} \nabla' \frac{\nabla m_\pi(Z; v)}{m_\pi(Z; v)}$$

(24)
$$= E_{\mu,v} \frac{1}{2} \left(\frac{\nabla^2 m_\pi(Z; v)}{m_\pi(Z; v)} - \|\nabla \log m_\pi(Z; v)\|^2\right).$$

Combining (20), (22) and (24) yields (18). That (18) equals (19) can be verified directly. □

It may be of independent interest to note that the intermediate step (21) is in fact a restatement of the well-known fact that any Gaussian convolution will solve the homogeneous heat equation, which has a long history in science and engineering; for example, see [20]. Brown, DasGupta, Haff and Strawderman [5] recently used identities derived from the heat equation, including one bearing a formal similarity to (21), in other contexts of inference and decision theory. Furthermore, as the Associate Editor kindly pointed out to us, the proof of Lemma 3 can also be obtained by appealing to Theorem 1 and equation (54) of that paper.

THEOREM 1. *Suppose $m_\pi(z; v_x)$ is finite for all $z$.*

(i) *If $\nabla^2 m_\pi(z; v) \leq 0$ for all $v_w \leq v \leq v_x$, then $p_\pi(y|x)$ is minimax under $R_{\mathrm{KL}}$. Furthermore, $p_\pi(y|x)$ dominates $p_{\mathrm{U}}(y|x)$ unless $\pi = \pi_{\mathrm{U}}$.*

(ii) *If $\nabla^2 \sqrt{m_\pi(z; v)} \leq 0$ for all $v_w \leq v \leq v_x$, then $p_\pi(y|x)$ is minimax under $R_{\mathrm{KL}}$. Furthermore, $p_\pi(y|x)$ dominates $p_{\mathrm{U}}(y|x)$ if for all $v_w \leq v \leq v_x$, $\nabla^2 \sqrt{m_\pi(z; v)} < 0$ on a set of positive Lebesgue measure.*

PROOF. As established in Corollary 1, $p_{\mathrm{U}}$ is minimax under $R_{\mathrm{KL}}$. Thus, minimaxity is established by showing that (16) is nonnegative, and dominance is establish by showing that (16) is strictly positive on a set of positive Lebesgue measure. Then (i) and (ii) follow from (18), (19) and the fact that $v_w < v_x$. □

COROLLARY 2. *If $m_\pi(z; v_x)$ is finite for all $z$, then $p_\pi(y|x)$ will be minimax if the prior density $\pi$ satisfies $\nabla^2 \pi(\mu) \leq 0$ a.e. Furthermore, $p_\pi(y|x)$ will dominate $p_{\mathrm{U}}(y|x)$ unless $\pi = \pi_{\mathrm{U}}$.*

PROOF. It is straightforward to show (see problem 1.7.16 of [15]) that $\nabla^2 m_\pi(z; v) \leq 0$ when $\nabla^2 \pi(\mu) \leq 0$ a.e. Therefore, Corollary 2 follows immediately from (i) of Theorem 1. □

The above sufficient conditions for minimaxity and domination in the KL risk prediction problem are essentially the same as those for minimaxity and domination in the quadratic risk estimation problem. What drives this connection is revealed by comparing Stein's unbiased estimate of quadratic risk reduction in (11) and (12) with (18) and (19). It follows directly from this comparison that the risk reduction in the quadratic risk estimation problem can be expressed in terms of $\log m_\pi$ as

$$(25) \qquad R_{\mathrm{Q}}(\mu, \hat{\mu}_{\mathrm{MLE}}) - R_{\mathrm{Q}}(\mu, \hat{\mu}_\pi) = -2 \left[ \frac{\partial}{\partial v} E_{\mu, v} \log m_\pi(Z; v) \right]_{v=1}.$$

**3. Examples.** In this section we show how Theorem 1 and Corollary 2 can be applied to establish the minimaxity of $\hat{p}_H$ and $\hat{p}_a$. Compared to the minimaxity proofs of Komaki [13] for $\hat{p}_H$, and of Liang [16] for $\hat{p}_a$, this unified approach is more direct and more general. We further indicate how our approach can be used to obtain wide classes of new minimax prediction densities.

EXAMPLE 1. Let us return to the Bayes predictive density $\hat{p}_H$, the special case of (3) under the harmonic prior $\pi_H(\mu)$ in (6). Following Komaki [13], the marginal of $Z|\mu \sim N_p(\mu, vI)$ under $\pi_H$ can be expressed as

$$(26) \qquad m_H(z;v) \propto v^{-(p-2)/2} \phi_p(\|z/\sqrt{v}\|),$$

where $\phi_p(u) = u^{-p+2} \int_0^{(1/2)u^2} t^{p/2-2} \exp(-t)\, dt$ is the incomplete Gamma function. By Lemma 2, $\hat{p}_H$ can be expressed in terms of this marginal as

$$(27) \qquad \hat{p}_H(y|x) = \frac{m_H(w; v_w)}{m_H(x; v_x)} \hat{p}_U(y|x).$$

Because $\pi_H$ is harmonic $[\nabla^2 \pi_H(\mu) \equiv 0$ a.e.$]$, and hence superharmonic, for $p \geq 3$, the fact that $\hat{p}_H$ is minimax and dominates $\hat{p}_U$ follows immediately from Corollary 2.

Beyond $\hat{p}_H$, one might consider the class of Bayes predictive densities $\hat{p}_\pi$ corresponding to the (improper) multivariate $t$ priors $\pi(\mu) = (\|\mu\|^2 + 2/a_2)^{-(a_1+p/2)}$. Because these priors are superharmonic for $a_1 \leq -1$ and $p \geq 3$, the minimaxity and domination of $\hat{p}_U$ by these rules follows immediately from Corollary 2.

EXAMPLE 2. Turning next to $\hat{p}_a$, the marginal of $Z|\mu \sim N_p(\mu, vI)$ under the Strawderman prior $\pi_a$ in (7) can be expressed as

$$(28) \qquad \begin{aligned} m_a(z;v) \propto \int_0^\infty & \left\{ 2\pi v \left( \frac{v_0}{v} s + 1 \right) \right\}^{-p/2} \\ & \times \exp\left\{ -\frac{\|z/\sqrt{v}\|^2}{2((v_0/v)s+1)} \right\} (s+1)^{a-2}\, ds. \end{aligned}$$

Because $\pi_H$ is the special case of $\pi_a$ when $a = 2$, it follows that $m_H(z;v)$ is the special case of $m_a(z;v)$ when $a = 2$. As Fourdrinier, Strawderman and Wells [6] showed, the marginal for any proper prior cannot be superharmonic, so that Theorem 1(i) cannot hold for $\hat{p}_a$ when $a < 1$. However, Theorem 1(ii) does hold for such $\hat{p}_a$, because $\sqrt{m_a(z;v)}$ is superharmonic for $v \leq v_0$ when $p = 5$ and $a \in [0.5, 1)$ or $p \geq 6$ and $a \in [0, 1)$. This fact can be obtained using $h(s) \propto (1+s)^{a-2}$ in Theorem 2 below, which extends Theorem 1 of [6].

THEOREM 2. *For a nonnegative function $h(s)$ over $[0, \infty)$, consider the scale mixture prior*

$$\pi_h(\mu) = \int \pi(\mu|sv_0) h(s)\, ds, \tag{29}$$

*where $\pi(\mu|sv_0) = N_p(0, sv_0 I)$. For $Z|\mu \sim N_p(\mu, vI)$, let*

$$m_h(z; v) \propto \int_0^\infty \{2\pi v(s+1)\}^{-p/2} \exp\left\{ -\frac{\|z/\sqrt{v}\|^2}{2(s+1)} \right\} rh(rs)\, ds \tag{30}$$

*be the marginal distribution of $Z$ under $\pi_h(\mu)$, where $r = v/v_0$. Let $h$ be a positive function such that:*

(i) $-(s+1)h'(s)/h(s)$ *can be decomposed as $l_1(s) + l_2(s)$, where $l_1 \leq A$ is nondecreasing while $0 < l_2 \leq B$ with $\frac{1}{2}A + B \leq (p-2)/4$,*

(ii) $\lim_{s \to \infty} h(s)/(s+1)^{p/2} = 0$.

*Then $\sqrt{m_h(z; v)}$ in (30) is superharmonic for all $v \leq v_0$, and when $v_x \leq v_0$, the Bayes predictive density $\hat{p}_h(y|x)$ under $\pi_h(\mu)$ in (29) is minimax.*

PROOF. The proof of Theorem 1 in [6] shows that $\sqrt{m_h(z; v_0)}$ in (30) is superharmonic when $v_0 = 1$, and it is straightforward to show that this is true for general $v_0$. From this fact, $\sqrt{m_h(z; v)}$ will be superharmonic for all $v \leq v_0$ if $h_r(s) := rh(rs)$ satisfies (i) and (ii) when $r \in (0, 1]$.

First we show that $h_r$ satisfies (i). By the assumptions on $h$, we have $-(s+1)h'(s)/h(s)$ decomposed as $\tilde{l}_1(s) + \tilde{l}_2(s)$. Then

$$-(s+1)\frac{h'_r(s)}{h_r(s)} = -\frac{r(s+1)}{rs+1}(rs+1)\frac{h'(rs)}{h(rs)}$$
$$= \frac{r(s+1)}{rs+1}[\tilde{l}_1(s) + \tilde{l}_2(s)].$$

Choose $l_i$ to be $\tilde{l}_i$ multiplied by $r(s+1)/(rs+1)$. They can be checked to satisfy the conditions since the factor $(rs+r)/(rs+1)$ is a nondecreasing function of $s$ and less than or equal to 1 when $0 < r \leq 1$. To see that $h_r$ satisfies (ii), note that

$$\frac{h_r(s)}{(s+1)^{p/2}} = \frac{h(rs)}{(rs+1)^{p/2}} r \left( \frac{rs+1}{s+1} \right)^{p/2}$$

goes to zero when $s \to \infty$ since the first term goes to zero by the assumption on $h$.

Thus $\sqrt{m_h(z; v)}$ will be superharmonic for all $v \leq v_0$. When $v_x \leq v_0$, the minimaxity of $\hat{p}_h(y|x)$ then follows from (ii) of Theorem 1. □

Going far beyond these results, Theorem 2 can be used to obtain wide classes of proper priors that yield minimax Bayes predictive densities $\hat{p}_h$. Following the development in Section 4 of [6], such $\hat{p}_h$ can be obtained with particular classes of shifted inverted gamma priors and classes of generalized $t$-priors.

**4. Further extensions.** Priors such as $\pi_H$ and $\pi_a$ are concentrated around 0, so that the risk reduction offered by $\hat{p}_H$ and $\hat{p}_a$ will be most pronounced when $\mu$ is close to 0. However, such priors can be readily recentered around a different point to obtain predictive estimators that obtain risk reduction around the new point. Because the superharmonicity of $m_\pi$ and $\sqrt{m_\pi}$ will be unaffected under such recentering, the minimaxity and domination results of Theorems 1 and 2 will be maintained. Minimax shrinkage toward a subspace can be similarly obtained by recentering such priors around the projection of $\mu$ onto the subspace.

To vastly enlarge the region of improved performance, one can go further and construct analogues of the minimax multiple shrinkage estimators of George [7, 8, 9] that adaptively shrink toward more than one point or subspace. Such estimators can be obtained using mixture priors that are convex combinations of recentered superharmonic priors at the desired targets. Because convex combinations of superharmonic functions are superharmonic, Corollary 2 shows that such priors will lead to minimax multiple shrinkage predictive estimators.

**Acknowledgments.** We would like to thank Andrew Barron, Larry Brown, John Hartigan, Bill Strawderman, Cun-Hui Zhang, an Associate Editor and anonymous referees for their many generous insights and suggestions.

E. I. GEORGE
DEPARTMENT OF STATISTICS
THE WHARTON SCHOOL
UNIVERSITY OF PENNSYLVANIA
PHILADELPHIA, PENNSYLVANIA 19104-6340
USA
E-MAIL: edgeorge@wharton.upenn.edu

F. LIANG
INSTITUTE OF STATISTICS
   AND DECISION SCIENCES
DUKE UNIVERSITY
DURHAM, NORTH CAROLINA 27708-0251
USA
E-MAIL: feng@isds.duke.edu

X. XU
DEPARTMENT OF STATISTICS
OHIO STATE UNIVERSITY
COLUMBUS, OHIO 43210-1247
USA
E-MAIL: xinyi@stat.ohio-state.edu